\documentclass[10pt]{article}                                       

\usepackage{amsfonts}

\title{On the Speed of the One-dimensional Excited Random Walk in the Transient Regime}
\author{T. Mountford, L.P.R. Pimentel, G. Valle}
\newtheorem{theorem}{Theorem}[section]                                          
\newtheorem{proposition}[theorem]{Proposition}                          
\newtheorem{lemma}[theorem]{Lemma}
\newtheorem{corollary}[theorem]{Corollary}

\newtheorem{remark}{Remark}

\newcommand{\no}{\noindent}
\newcommand{\nn}{\nonumber}
\newcommand{\oo}{\infty}

\newcommand{\ra}{\rightarrow}

\newcommand{\Z}{\mathbb{Z}}

\newcommand{\N}{\mathbb{N}}

\newcommand {\F}{\mbox{${ F}$}}

\begin{document}

\maketitle

\begin{abstract}
We study a class of nearest-neighbor discrete time integer random walks introduced by Zerner, the so called multi-excited random walks. The jump probabilities for such random walker have a drift to the right whose intensity depends on a random or non-random environment that also evolves in time according to the last visited site. A complete description of the recurrence and transience phases was given by Zerner under fairly general assumptions for the environment. We contribute in this paper with some results that allows us to point out if the random walker speed is strictly positive or not in the transient case for a class of non-random environments.
\end{abstract}

%%%%%%%%%%%%%%%%%%%%%%%%%%%%%%%%%%%%%%%%%%%%%%%%%%%%%%%%%%%%%%%
\section{Introduction} 
%%%%%%%%%%%%%%%%%%%%%%%%%%%%%%%%%%%%%%%%%%%%%%%%%%%%%%%%%%%%%%%
\setcounter{equation}{0}

The excited (or multi-excited if we follow Zerner's \cite{zerner} terminology) random walk (ERW) on $\Z^d$ can be informally described as a model of a random walk on $\Z^d$ that upon its first $M_x$ visits to site $x$ the walk is pushed toward a specific direction and on subsequent visits a neighbor is chosen uniformly at random. This model was introduced by \cite{benjamini-wilson} and it gives an interesting intermediary class between classical random walks and random walks in random environment. For the sake of simplicity, we adopt Zerner's terminology \cite{zerner} by saying $M_x$ is the initial number of cookies in site $x$, and that upon eating a cookie, i.e., visiting a site which still has cookies, the random walk has a positive drift, say to the "right". The determination of the initial number of cookies and the intensity of the drift imposed on the walk each time a cookie is eaten, which can be random or non-random, gives what we call the environment for the excited random walk. 

The first non trivial question about the behavior of the ERW is if it is transient. This was first studied by Benjamini and Wilson \cite{benjamini-wilson}, where it was proved that the non-random ERW is transient in dimension $d>1$ if there is at least one cookie per site. The further question is naturally about the strict positivity of the walk's speed in the transient regime, for one or more cookies in dimension 3 this was affirmatively answered by Kozma in \cite{kozma}. While in one dimension, for at most one cookie per site it is simple to show recurrence, it is not evident the behavior of the walk with two or more cookies per site. Transient and recurrent regimes in one-dimension were completely characterized by Zerner \cite{zerner} not only for non-random ERW, but also for those with an stationary and ergodic environment. In his paper Zerner also shows that the one-dimensional ERW with two or less cookies per site has null speed, however the same question remains open for three or more cookies per site. Our aim is to give a partial answer to this question. 

To present our results, we start with the formal description of the ERW. An environment is of the form $\omega = (\omega(x))_{x\in \Z}$, where for each $x \in \Z$, $\omega \in [1/2,1]^\N$. The discrete time ERW with starting point and environment $(x_0,\omega)$ is an stochastic process $(X_n,\omega_n)$ whose transition is defined by the following rule
$$
\left\{
\begin{array}{ll}
(X_n,\omega_n) \ra (X_n - 1,\theta_{X_n} \omega_n) & , \textrm{ with probability } 1 - \omega_n(X_n,1) \\
(X_n,\omega_n) \ra (X_n + 1,\theta_{X_n} \omega_n) & , \textrm{ with probability } \omega_n(X_n,1) \, .
\end{array}
\right.
$$ 
where $\theta_z: [1/2,1]^{\Z \times \N} \ra [1/2,1]^{\Z \times \N}$ is defined by
$$
\theta_z \omega (x,n) = \left\{
\begin{array}{ll}
\omega(x,n) & , \ x \neq z \\
\omega(z,n+1) & , x=z
\end{array}
\right.
$$
for every $n \in \N$.

\medskip

We will make two further assumptions on an environment: First that on each site there is initially only a finite number of cookies and second that on each visit to a site with available cookies, a cookie must be consumed. Therefore for such an environment $\omega$ we have that there exists an $M_x \in \N$ such that $\omega(x,n) > 1/2$ if and only if $n\le M_x$. For a $p>1/2$ and $M \in \N$, we will denote by $\omega^{M,p}$ the homogeneous environment with $M_x = M$, for every $x$, and $\omega(x,n) = p$ for $n\le M$. 

Following \cite{zerner}, an environment $\omega$ as above will be called a cookie environment. This is reasonable from the point of view that if $\omega(x,n) > 1/2$ we say that a cookie of intensity $\omega(x,n)$ will be ready to be consumed by the time of the $n$-th visit to the site $x$. 

In \cite{zerner} it was shown that the excited random walk starting at $(0,\omega^{M,p})$ is transient if and only if $M(2p-1)>1$ and he obtained a law of large numbers for the excited random walk under more general environments proving that the walk's speed
$$
\textrm{a.s-}\lim_{n \ra \oo} \frac{X_n}{n}
$$
is well defined. We also find in \cite{zerner} that for $p<1$ and $M=2$ the walk's speed is $0$. The question of what the walk's speed is, in the transient regime for homogeneous environments with more than 2 cookies per site remains open. Our results
in this direction are the following:

\medskip

\begin{theorem}
\label{thm1} 
For the ERW starting at $(0,\omega^{M,p})$ we have that
\begin{enumerate}
\item[(i)] For every $p\in (1/2,1)$, there exists $M_0 = M_0(p)$ sufficiently large such that the walk's speed is strictly positive for all $M>M_0$.
\item[(ii)]  If the bias $p$ and the number of cookies, $M$, satisfy
$$
M(2 p-1) \in (1,2)
$$
then the cookie r.w. is transient but with speed 0.
\end{enumerate}
\end{theorem}

\medskip

\begin{remark}
A counter-example is given in section \ref{section:thm1(i)} where for all $M>0$ condition (i) does not holds for an ergodic environments with $M$ as the mean number of cookies per site.
\end{remark}

We do believe that the result in (ii) of Theorem \ref{thm1} is almost optimal, i.e., that if $M(2p-1) > 2$ then the ERW starting at $(0,\omega^{M,p})$ has positive speed. In this direction we show that:

\medskip

\begin{theorem}
\label{thm2}
There exist constants $C_0>1$ and $0<\epsilon_0 < \frac{1}{2}$ such that if $p \le \frac{1}{2} + \epsilon_0$ and $M(2p-1) \geq C_0$ then the ERW starting at $(0,\omega^{M,p})$ has strictly positive speed.
\end{theorem}

%%%%%%%%%%%%%%%%%%%%%%%%%%%%%%%%%%%%%%%%%%%%%%%%%%%%%%%%%%%%%%%
\section{Technical Estimates} 
\label{section:techest}
%%%%%%%%%%%%%%%%%%%%%%%%%%%%%%%%%%%%%%%%%%%%%%%%%%%%%%%%%%%%%%%
\setcounter{equation}{0}

We start this section with some definitions, which will not be used only here but through the entire paper. Consider an ERW $(X_n,\omega_n)_{n\in \N}$ starting at some $(0,\omega)$. For integers $R,S>0$ let $T_R = \inf \{ n \ge 1 : X_n = R\}$, $T_S^R = \inf \{ n \ge T_R : X_n = S\}$ and $N_R$ ($N_R^n$) be respectively the number of visits of the walk to site $R$ (before time $n$). For an environment $\omega$ and $x \in \Z$ we write $\mathrm{P}_{(x,\omega)}$ for the probability induced by the excited random walk starting at $(x,\omega)$ and $E_{(x,\omega)}$ for the corresponding expectation operator.  When $\omega$ is the environment with no cookies the exited random walk is a simple random walk and we only write $P_x$. Finally let $(\mathcal{F}_n)_{n\ge 1}$ be the filtration generated by $(X_n)_{n\ge 0}$.

\medskip

Denote by $D_k$ the total drift of the ERW $(X_n,\omega_n)_{n\ge 1}$ accumulated until time $k$, i.e., $(2p-1)$ times the number of cookies eaten before time $k$. Let us recall another observation of \cite{zerner} that $V_n = X_n - D_n$ is a $\mathrm{P}_{(0,\omega)}-$martingale with respect to the filtration generated by $(X_n)_{n\ge 1}$. Therefore we obtain the martingale decomposition for $X_n$ as $V_n + D_n$, where $V$ is a martingale whose increments are either $\pm 1$ or $2(1-p)$ or $-2p$ and $D$ is an increasing process, with both processes starting at value 0.

\medskip

\begin{lemma} \label{lem1} For an ERW starting at the origin with $M$ cookies in at least $\alpha N$ sites in $(-N,0]$, the probability that the cookie r.w. reaches $-N$ before eating $\frac{M{\alpha} N}{4}$ cookies is bounded by $\frac{2M}{\alpha N}$
\end{lemma}

\medskip

\no \textbf{Proof:}
There are at least $\frac{\alpha N}{2}$ sites in $(-N(\frac{1-\alpha}{2}), 0)$ with $M$ cookies. We denote these sites from right to left by $x_1, x_2, \ldots, x_R$  for $R = \frac{\alpha N}{2}$. For $x_i$ fixed and based on monotonicity results obtained in \cite{zerner}, we consider without loss of generality the cookie r.w. $X$ with no cookies to the right of $x_i$. At each visit to $x_i$ the probability that the cookie r.w. returns to $x_i$ before hitting $-N$ is greater or equal to the same probability computed with respect to a simple symmetric random walk which is given by $\frac{1}{2} + \frac{1}{2} \left( 1 - \frac{2}{\alpha N} \right) = 1- \frac{1}{\alpha N}$. Thus the probability that $M$ cookies are eaten at $x_i$ before hitting $-N$ is $(1- \frac{1}{\alpha N})^M$. So if $V$ is the number of $x_i$ at which less than $M$ cookies are eaten before hitting $-N$. Then
$$
E(V) \leq \frac{\alpha N}{2} \left(1-\left(1- \frac{1}{\alpha N}\right)^{\!\!\! M} \right),
$$
so $\mathrm{P}(V \geq \frac{\alpha N}{4}) \leq 2(1-(1- \frac{1}{\alpha N})^M) \leq \frac{2M}{\alpha N}$. $\square$

\bigskip

\begin{lemma} 
\label{lem2}
For an ERW starting at the origin with $M$ cookies initially present for at least $\alpha N$ sites in $(-N, 0]$, the probability that the cookie r.w. reaches $-N$ before $N$ is bounded by
$$
\frac{2 M}{\alpha N} + \exp\left( - \, c \, \frac{M\alpha}{4}(2p-1) \right)
$$
for $c$ not depending on $p$ or $N$ or $M$ or $\alpha$.
\end{lemma}

\medskip

\no \textbf{Proof:}
Consider the martingale decomposition for the position of the ERW $X_n = V_n + D_n$ as above. Now we have the inclusion $\{ T_{-N} < T_N\}$ is contained in the union $\{X$ `` eats" $ < \frac{M{\alpha} N}{4}$ cookies before hitting $-N\} \cup \{ \textrm{There exists no } n < n^\prime \leq S_{N,V} : V_{n^\prime}-V_n \geq 2N\}$, where $S_{N,V}$ is the first time $V$ touches $-\left( \frac{M \alpha}{4} +1 \right)(2p-1)N$. By standard Brownian embedding we have that the probability of the later event is bounded by
$$
e^{-c(\frac {M{\alpha}}{4} + 1)(2p-1)}
$$
for some $c$ universal not depending on $p$ or $M$ and the result follows.
$\Box$  

\bigskip

\begin{corollary} \label{cor1}
For every $\epsilon > 0$, there exists $M_1$ and $n_1$ so that for all $N  \geq n_1$, if the ERW begins at zero with an environment $\omega$ such that $\frac14 $ of the sites in $(-N,0)$ have at least $M_1$ cookies in them, then
$$
\mathrm{P}_{(0,\omega)} (T_{N} < T_{-N}) \ \geq \ 1- \epsilon .
$$
\end{corollary}

\medskip

By repeated application of the Strong Markov property we then
obtain

\smallskip

\begin{corollary} \label{cor2}
For every $\epsilon > 0$, there exists $M_1$ and $n_1$ so that for all $N \ge n_1$ and positive integer $k_0$, if the ERW begins at zero with an environment $\omega$ such that for $k = 0,1,2 \cdots k_0$, the number of sites with at least $M_1$ cookies in spatial interval $((-2^{k+1}+1)N,(-2^k+1)N]$ exceeds $ 2^{k-1} N $, then
$$
\mathrm{P}_{(0,\omega)}(\, T_{-2^{k_0}N} \le T_N \, ) \ \leq \  \epsilon^{k_0} .
$$
\end{corollary}

\medskip

The following is crude but useful.
\begin{lemma} 
\label{lem3} 
Let $(X_n : n \geq 0)$ be an ERW with $X_0 = 0$ and initial environment so that in interval $(-c N, N)$ there are less than $\gamma N$ cookies. The probability that $X$ hits $-cN$ before $N$ is at least
$$
\frac{1-\gamma(2p-1)}{1+ c + 2/N}.
$$
\end{lemma}

\smallskip

\no \textbf{Proof:} Until leaving interval $(-cN, N)$ the position of the ERW $X_n$ is equal to $V_n + D_n$ (see the beginning of this section).  Thus the event $\{T_{-cN} \le T_N \}$ contains the event that $V \textrm{ hits } (-\infty, - (c+\gamma(2p-1))N] \textrm{ before } [N(1-\gamma(2p-1)), \infty)$. Now for $\tau$ the first hitting time for $V$ of $(-\infty, - (c+\gamma(2p-1))N] \cup [N(1-\gamma(2p-1)), \infty)$ we have
\begin{eqnarray*}
0 = E(V_\tau) &\geq& (-(c+\gamma(2p-1))N-2)\mathrm{P}(V_\tau \leq -(c+\gamma(2p-1))N)\\
&+& (N(1-\gamma(2p-1))(1-\mathrm{P}(V_\tau \leq -(c+\gamma(2p-1))N))
\end{eqnarray*}
so
$$
\mathrm{P}(V_\tau \leq - (c+\gamma (2p-1))N) \geq \frac{N(1-\gamma(2p-1))}{N(1+c)+2} \qquad \Box
$$

\medskip

We similarly have

\smallskip

\begin{lemma} 
\label{lem4}
Consider an ERW $(X_n,\omega_n)_{n \geq 0}$  so that for the initial environment $\omega$ there are less than $\gamma N$ cookies on $(-N,N)$, for $\gamma(2p-1) < 1$. Then the probability that $X$ exits $(-N, N]$ after time $N^2/2$ is at least $c_0$ for some constant $c_0$ depending only on $\gamma$.
\end{lemma} 

\medskip

We now seek to refine Lemma \ref{lem3}.

\begin{lemma} 
\label{lem5} 
Fix $\epsilon > 0$ and $N < \infty$ such that $\epsilon (2p-1) \le N^{-2}$. Let $\Omega$ be the event: the ERW starting at 0
\begin{enumerate}
\item[(a)] exits $(-2^k , 2^k)$ to the left and 
\item[(b)] on exiting this interval, the number of cookies in $[-2^k, 2^k)$ is less than $\epsilon \, 2^{k+1}$.
\end{enumerate}

Define $C_k = C_k(\epsilon, N, b,M)$ as the infimum of $\mathrm{P}_{(0,\omega)}(\Omega)$ over all environments $\omega$ having at most $M$ cookies per site and satisfying
\begin{itemize}
\item[(i)] For every $0 \leq i < N$ the number of cookies in $(-\frac{(i+1)2^k}{N} , \frac{-i 2^k}{N}]$ is less or equal to $\frac{b2^k}{N}$.
\item[(ii)] The number of cookies in $[0,2^k)$ is less or equal to $\epsilon 2^k$.
\end{itemize}

Then $\lim_{k \to \infty} \, C_k \geq  2^{-(1+b(2p-1))} (1+c_N)$, where $c_N \ra 0$ as $N \ra +\oo$.
\end{lemma}

\smallskip

\no \textbf{Proof:} Fix an environment $\omega$ satisfying (i) and (ii) in the statement. We consider the event $A(i) = \{ T_{- \frac{i 2^k}{N}} \le T_{2^k} \}$ and let $B(i)$ be the event that on hitting  $-\frac{i2^k}{N}$ there are less than $\frac{\epsilon}{2} \frac{2^k}{N}$ cookies remaining in interval $\left(- \frac{i 2^k}{N}, -\frac{(i-1) 2^k}{N} \right]$. By the strong Markov property we can write
\begin{eqnarray*}
\mathrm{P}_{(0,\omega)}(\Omega) \ge \prod_{i=1}^{N} \mathrm{P}_{\left(-\frac{(i-1)2^k}{N},\omega_{\tilde{T}_{i-1}} \right)} \left(A(i) \, | \cap^{i-1}_{j=1} B(j)\right) - \sum_{i=1}^{N} \mathrm{P}_{(0,\omega)} (A(i)/B(i)),
\end{eqnarray*}
where $\tilde{T}_{i-1}$ is equal to $T_{-\frac{(i-1)2^k}{N}}$. 

\medskip

We first consider the event $A(i)/B(i)$. To bound this probability we follow along the lines of Lemma \ref{lem3}. We consider for each $x \in I_k = \left[\frac{-i2^k}{N} + 2^{k/2}, -\frac{(i-1)2^k}{N}\right)$ the event $H_x$ that all cookies at $x$ are eaten before $X$ hits $-\frac{i2^k}{N}$. Now observe that $P(H_x) \to 1$ for $k \to \infty$, uniformly in $x$ and environment. Thus the expected number of cookies remaining in $\left(-\frac{i2^k}{N}, - \frac{(i-1)2^k}{N}\right]$ when $X$ hits $-\frac{i2^k}{N}$ is bounded by
$$
M 2^{k/2} + \frac{2^k}{N} \left(\sup_{x\in I_k} \mathrm{P}(H^c_x)\right) = o(2^k).
$$
Thus $\sup_i \mathrm{P}(A(i)/B(i)) \to 0$ as $k \to \infty$ uniformly over environments satisfying (i) and (ii).

\medskip

We now obtain bounds on $\mathrm{P}\left(A(i) \left\vert \cap^{i-1}_{j=1} B(j) \right. \right)$
for $i=1$. We apply Lemma \ref{lem3} with 
$$
N = 2^k \left( 1+ \frac{i-1}{N} \right) \, , \qquad c= \frac{1}{N+i-1} \, \quad \textrm{ and } \quad \gamma = \frac{\epsilon + \epsilon/2 \frac{i-1}{N} + \frac{b}{N}}{1 + \frac{i-1}{N}} <  \epsilon + \frac{b}{N+i-1}.
$$
This give a lower bound of
\begin{eqnarray*}
\lefteqn{\!\!\!\!\!\!\!\!\!\!\!\!\!\!
\frac{1- \epsilon(2p-1) - \frac{b(2p-1)}{N+i-1}}{\frac{N+i}{N + i-1}}
= \frac{(N+i)(1- \epsilon(2p-1)) - b(2p-1)}{N+i}}\\
&=& 1 - \frac{(1+b(2p-1))}{N+i} - \frac{N+ i-1}{N+i} \epsilon(2p-1)\\
&\geq& 1 - \frac{(1+ b(2p-1))}{N+i} - \frac{\epsilon(2p-1)}{N}.
\end{eqnarray*}

From this we obtain that
\begin{eqnarray*}
\lim_{k\to \infty }C_k &\geq& \prod^N_{i=1} \left(1- \frac{1+b(2p-1))}{N+i} - \frac{\epsilon (2p-1)}{N} \right)\\
&=& 2^{-(1+b(2p-1))} (1-\mathbf{o}(1)) \, . \quad \Box
\end{eqnarray*}

\bigskip \bigskip

%%%%%%%%%%%%%%%%%%%%%%%%%%%%%%%%%%%%%%%%%%%%%%%%%%%%%%%%%%%%%%%
\section{Proof of (i) in Theorem \ref{thm1}} 
\label{section:thm1(i)}
%%%%%%%%%%%%%%%%%%%%%%%%%%%%%%%%%%%%%%%%%%%%%%%%%%%%%%%%%%%%%%%
\setcounter{equation}{0}

Consider the ERW $(X_n,\omega_n)_{n\in \N}$ starting at $(0,\omega^{M,p})$ for some $M>0$ and $0<p<1$.
As a consequence of Zerner's result \cite{zerner} on the existence of the ERW speed we also have that $\frac{T_R}{R}$ converges almost surely to the inverse of the walk's speed. Since
$$
T_R = \sum_{x\le R} N_x^{T_R} \, ,
$$
to obtain the desired result, by Fatou's lemma, we only have to show that there exists a sequence $(R_i)_{i=1}^{+\oo}$ with $R_i \ra \oo$ such that there exists a sequence $C>0$ with
$$
\sup_i \frac{1}{R_i} \sum_{x\le R_i} E[ N_x^{T_{R_i}} ] \le C \, .
$$
Our aim will be to show the stronger result that
$$
\sup_x E [ N_x ] < + \oo \, ,
$$
which follows from the following proposition:

\medskip

\begin{proposition}
\label{prop:estNx}
For every $p \in (1/2,1)$ there exists a $M_0=M_0(p)$ sufficiently large such that for 
the excited random walk $(X_n,\omega_n)_{n\in \N}$ starting at $(0,\omega^{M,p})$ with $M>M_0$ we have for every $n>1$ 
$$
\sup_{x \ge 1} \mathrm{P} ( N_x \ge n ) \le \frac{C}{n^\gamma} 
$$
for some constants $C>0$ and $\gamma > 1$.
\end{proposition}

\medskip

Our first step to estimate the supremum in the statement of the previous proposition is to use coupling arguments and choose an
appropriate environment $\tilde{\omega}^{M,p}$ such that 
$$
\mathrm{P}_{(0,\omega^{M,p})} ( N_x \ge n ) \le \mathrm{P}_{(0,\tilde{\omega}^{M,p})} ( N_0 \ge n )
$$
for every $x \in \N$. We take as the environment $\tilde{\omega}^{M,p}$ the one with no cookies strictly to the left of the origin and $M$ cookies of intensity $p$ to the right of the origin. Clearly this choice satisfies the previous inequality. 

\medskip

\begin{lemma}
\label{lemma:contvis}
For every $p\in (1/2,1)$, there exists $M_0 = M_0(p)$ sufficiently large such that for all $M>M_0$ for the excited random walk $(X_n,\omega_n)_{n\in \N}$ starting at $(0,\tilde{\omega}^{M,p})$ we have for every $n>1$ 
$$
\mathrm{P}_{(0,\tilde{\omega}^{M,p})} ( N_0 \ge n ) \le \frac{C}{n^\gamma} 
$$
for some constants $C>0$ and $\gamma > 1$.
\end{lemma}

\medskip

For an arbitrary $R>0$ we have that the probability in the statement of the previous lemma is bounded above by
\begin{equation}
\label{eq:contvis1}
\mathrm{P}_{(0,\tilde{\omega}^{M,p})} ( N_0^{T_R} \ge n ) +
\mathrm{P}_{(0,\tilde{\omega}^{M,p})} ( T_0^R < +\oo ) \, .
\end{equation}
The first term in the previous expression is dominated by (recall $P_x$ refers to simple random walk probabilities started at $x$)
$$
\mathrm{P}_0 ( N_0^{T_R} \ge n ) = \mathrm{P}_0 ( T_0 \le T_R )^n \le \left( 1 - \frac{c}{R} \right)^n 
$$
for some constant $c>0$ not depending on $R$, this follows by comparison with the simple symmetric random walk, see Lemma 1 in \cite{zerner}. From this moment on we take $R=\left\lfloor n^\alpha \right\rfloor$ for some $0 < \alpha < 1$, from where a bound as the one in the statement of Lemma \ref{lemma:contvis} is easily obtained for the first term in (\ref{eq:contvis1}). 

We will have more work to deal with the second term in (\ref{eq:contvis1}) and obtain the bound
\begin{equation}
\label{contvis2}
\mathrm{P}_{(0,\tilde{\omega}^{M,p})} ( T_0^{\left\lfloor n^\alpha \right\rfloor} < +\oo ) \le \frac{C}{n^\gamma}
\end{equation}
for some $C>0$ and $\gamma > 1$. 

\medskip

To begin the proof we fix some events: Take $\epsilon > 0$ small, which will be fixed later, and choose $M_1 = M_1(\epsilon)$ large enough in the sense of Corollary \ref{cor1}. We define the event $\Gamma(e, n)$ that the ERW starting at $(0,\tilde{w}^{M,p})$ satisfies
\begin{enumerate}
\item[(a)] after hitting $2n$ for the first time, it hits $4n$ before hitting $n$ and
\item[(b)] upon hitting $4n$, there are more than $e \cdot n$ sites in interval $(n, 2n)$ with at least $M_1$ cookies. 
\end{enumerate}

For $x \geq 0$, we denote $\Gamma^x (e, n)$ the ``shift" of event $\Gamma(e, n)$ that is the intersection of
\begin{enumerate}
\item[(a)] after time $T_{x+2n} \, (= T^x_{x+2n})$, the ERW hits $x+4n$ before $x+n$
\item[(b)] at time $T_{x+4n} \, (=T^x_{x+4n})$ the number of sites in $(x+n, x+2n)$ with at least $M_1$ cookies is at least $e \cdot n$. 
\end{enumerate}

As a consequence of Lemma 10 in \cite{zerner}, note that   
$$
\mathrm{P}_{(0,\omega)} (\Gamma^x(e, n)| \mathcal{F}_{T_x}) = \mathrm{P}_{(0,\omega)}(\Gamma(e, n)).
$$
for all environments $\omega$ with fully occupied (M cookies per site) to the right of $x$ (included). Note also that $\Gamma^x(e, n)$ and $\Gamma^y(e, n)$ are independent for $y < 0 \leq x \leq y -4n$ under $\mathrm{P}_{(0,\omega)}$.

\medskip

Now let $(n_i)_{i=1}^{+\oo}$ be an increasing sequence of integers such that $n_1 = n_1(\epsilon)$ is fixed large enough to satisfy Corollary \ref{cor1}, but for the moment not fully specified, and $n_{j+1} = \left\lfloor n_j^{3/2} \right\rfloor$ for $j\ge 1$. Observe that for an arbitrary sequence $(e_i)_{i=0}^{+\oo}$, we have, provided $n_1$ is sufficiently large, that
$$
\bigcap_{j=h(n^\alpha)}^{+\oo} \bigcap_{l=0}^{3 \left\lfloor n_j^{1/2} \right\rfloor} \Gamma^{l n_j} (e_j,n_j) \subset \{ T_0^{\left\lfloor n^\alpha \right\rfloor} = +\oo \} \, ,
$$
where $h(n^\alpha) = \sup\{j:n_j < \frac{{\left\lfloor n^\alpha \right\rfloor}}{2}\}$. Therefore
\begin{eqnarray}
\label{eq:contvis3}
\lefteqn{ \!\!\!\!\!\!\!\!\!\!\!\!\!\!\!\!\!\!\!\!\!\!\!\!\!\!\!\!\!\!\!\!\!\!\!\!\!\!\!
\mathrm{P}_{(0,\tilde{\omega}^{M,p})} ( T_0^{\left\lfloor n^\alpha \right\rfloor} < +\oo ) \le
\sum_{j=h(n^\alpha)}^{+\oo} \sum_{l=0}^{3 \left\lfloor n_j^{1/2} \right\rfloor} \mathrm{P}_{(0,\tilde{\omega}^{M,p})}  \left( \Gamma^{ln_j} (e_j,n_j)^c \right)} \nn \\
& & \le 
\sum_{j=h(n^\alpha)}^{+\oo} (3 n_j^{1/2}+1) \, \mathrm{P}_{(0,\tilde{\omega}^{M,p})}  \left( \Gamma (e_j,n_j)^c \right) 
\end{eqnarray}
where the second inequality follows from monotonicity arguments (see Section 7 in \cite{zerner}).

\medskip

We shall prove the following result:

\medskip 

\begin{lemma}
\label{lemma:gamest} 
There is a
 decreasing sequence $(e_i)_{i=1}^\oo$ with $e_1 = \frac34 , \ e_i \ge 1/2$ for every $i$ and $n_j$ as above such that there exists a $M_0$ sufficiently large so that for all $M>M_0$ and every $j\ge 1$,
\begin{equation}
\label{eq:ind}
\mathrm{P}_{(0,\tilde{\omega}^{M,p})}  \left( \Gamma (e_j , n_j)^c \right) \le \frac{1}{n_j^\zeta}  
\end{equation}
for some $\zeta > \frac{3 + \alpha}{2\alpha}$.
\end{lemma}

\medskip

From Lemma \ref{lemma:gamest} and (\ref{eq:contvis3}), we obtain (\ref{contvis2}). We are going to show Lemma \ref{lemma:gamest} by induction and the first step is:

\medskip

\begin{lemma} 
\label{lemma:gamest1}
For every $1/2 < p < 1$ and $\zeta > 1$, we have that for every $R$ sufficiently large there exists $M_0 = M_0(n, p, \zeta)$ such that for $M > M_0$
$$
\mathrm{P}_{(0,\tilde\omega^{M,p})} (\Gamma(3/4, R)^c) \leq \frac{1}{R^\zeta} \, .
$$
\end{lemma} 

\smallskip

\no \textbf{Proof:} We couple the ERW $(X_n)_{n=1}^{+\oo}$ with the symmetric random walk reflected at the origin $(Y_m)^{+\infty}_{n=1}$ as in Lemma 1 in \cite{zerner}  . We have almost surely that $X_n \ge Y_n$, $n \geq 1$. Therefore,
for $M_0 > M_1=M_1(\epsilon)$
\begin{eqnarray*}
\mathrm{P}_{(0, \tilde \omega^{M,p})} (\Gamma(3/4, R)^c) &\leq&  \mathrm{P}_{(0, \tilde \omega^{M,p}),0} \left(\Gamma(3/4, R)^c \cap \left\{ \sup_{n \leq (M_0-M_1)} Y_n < 4R \right\}\right) +\\
&& + \mathrm{P}_{(0, \tilde \omega^{M,p}),0} \left(\Gamma(3/4, R)^c \cap \left\{\sup_{n \leq (M_0-M_1)} Y_n \geq 4R\right\}\right)\\
&& \leq \mathrm{P}_0 \left(\sup_{n \leq (M_0-M_1)} Y_n < 4R \right)+\\
&& + \mathrm{P}_{(0, \tilde\omega^{M,p}),0} \left(\Gamma(3/4, R)^c \left\vert \left\{\sup_{n\leq (M_0-M_1)} Y_n \geq 4R \right\} \right. \right).
\end{eqnarray*}

Now suppose $M_0 \ge R^{2(2\zeta^\prime+1)} + M_1$ then by the CLT the first probability in the rightmost term of the previous expression is of order $R^{-2\zeta}$. Thus it remains to deal with the second probability. Now observe that if $\{\sup_{n \le (M_0 - M_1)} Y_n \geq 4R\}$ happens then the excited random walk has arrived at $4R$ without eat all the cookies on a single site in the interval $[0,4R]$ and therefore the second probability is bounded above by the probability that an asymmetric simple random walk with jump probability $p$ starting at $2R$ touches $R$ before $4R$, this is the Gambler's ruin probability which is bounded above by $e^{-cR}$ for some constant $c$ depending on $p$ but not on $R$. Now taking $R$ sufficiently large we obtain the inequality in the statement. $\square$

\bigskip

\no \textbf{Proof of Lemma \ref{lemma:gamest}:} Apply Lemma \ref{lemma:gamest1} with $R=n_1$ to obtain an $M_0=M_0(n_1, p, \zeta)$ such that for all $M>M_0$, (\ref{eq:ind}) holds for $j=1$. Then suppose that (\ref{eq:ind}) holds for some arbitrary integer $i \ge 1$, we will show that it also holds for $i+1$. To guarantee the uniform bound, we should later adjust $n_1$ by making it larger if necessary.

\smallskip

We first note that if the event $\Gamma^x(e_i, n_i)$ occurs this says very little, in principle, about the number of cookies in $(x+n_i, x+2n_i)$ by the time the point $y$ is reached for $y > x + 4n_i$. However if $\Gamma^x(e_i, n_i) \cap \Gamma^{x+n_i}(e_i, n_i)$ occurs then after hitting $x+3n_i$ for the first time (necessarily between $T_{x+2n_i}$ and $T_{x+4n_i}$) the cookie r.w. must hit $x+5n_i$ before hitting $x+2n_i,$. If
$$
\Gamma^x(e_i, n_i) \cap \Gamma^{x+n_i}(e_i, n_i) \cap \Gamma^{x+2n_i}(e_i, n_i) 
$$
occurs then after hitting $x+4n_i$ the cookie r.w. must hit $x+6n_i$ before $x+3n_i$ and hence before $x+2n_i$. By induction we obtain that on 
$$
\Gamma^x(e_i, n_i) \cap \Gamma^{x+n_i}(e_i, n_i) \ldots \cap \Gamma^{x+rn_i}(e_i, n_i)
$$
the ERW hits $x+(r+4)n_i$ before $x+n_i$ after hitting $x+2n_i$ for the first time and at time $T_{x+(r+4)n_i}$ the number of sites in $(n_i, 2n_i)$ with $M_0$ or more cookies is at least $e_i n_i$. Since for every $0 < s < r$ we have
$$
\bigcap^r_{j=0} \Gamma^{x+jn_i}(e_i,n_i) \subseteq \bigcap_{j=s}^r \Gamma^{x+jn_i}(e_i, n_i)
$$
we have that on
$$
\bigcap^r_{j=0} \Gamma^{x+jn_i}(e_i, n_i),
$$
for every interval $(x+jn_i, x+(j+1)n_i), \  j = 1, \ldots , r$, the number of sites with more than $M_0$ cookies is at least $e_i n_i$ and 
$$
T_{x+(j+1)n_i}^{x+(j+2)n_i} > T_{x+(4+r)n_i}^{x+(j+2)n_i}.
$$

\medskip

Now based on the previous discussion we also have that if $e_{i+1} \leq e_i$ than, if $n_1$ is large enough,
$$
\Gamma(e_{i+1}, n_{i+1}) \supset \Gamma(e_i, n_{i+1}) \supseteq \bigcap^{4 \frac{n_{i+1}}{n_i}}_{l= \frac{n_{i+1}}{2n_i}}  \Gamma^{ln_i} (e_i, n_i) \, .
$$
This implies that for $\Gamma(e_{i+1}, n_{i+1})^c$ to happen it is necessary that at least one of $\Gamma^{jn_i} (e_i,n_i)^c$ also happens for $\frac{n_{i+1}}{2n_i} \leq k \leq \frac{4n_{i+1}}{n_i}$. But, by independence, the probability that there exist $0 \leq k \leq j-4 \leq \frac{4n_{i+1}}{n_i} - 4$ such that $\Gamma^{kn_i}(e_i , n_i)^c \cap \Gamma^{jn_i} (e_i,n_i)^c$ occurs is bounded above by
\begin{equation}
\label{eq:twobad}
\left( \frac{4n_{i+1}}{n_i} \mathrm{P}(\Gamma(e_i, n_i)^c) \right)^2 .
\end{equation}
This now reduces the estimate on $\mathrm{P}(\Gamma(e_{i+1}, n_{i+1})^c)$ to dealing with the event where there are at most four consecutive $l \in [\frac{n_{i+1}}{2n_i} - 3, \frac{4n_{i+1}}{n_i}]$ with $\Gamma^{ln_i} (e_i, n_i)^c$ occurring. We will see that we need essentially to consider that $\Gamma^{ln_i} (e_i, n_i)^c$ occurs for a single $l$. We now consider the event $B_1$ that there exists a unique $l$ so that
\begin{enumerate}
\item[(a)] $l \in \left(\frac{n_{i+1}}{2n_i} - 3, 4 \frac{n_{i+1}}{n_i} \right)$
\item[(b)] $\Gamma^{ln_i} (e_i, n_i)^c$ occurs
\item[(c)] $\Gamma^{jn_i} (e_i, n_i)$ occurs for every $j \in \left( l+3, 4\frac{n_{i+1}}{n_i} \right)$ and for every $j \in (0, l).$
\end{enumerate}
On this event we note that there are two ways that $\Gamma^{kn_i}(e_i, n_i)$ cannot happen: if too many cookies in $((k+1)n_i, (k+2)n_i)$ are eaten, which necessarily occurs while the ERW $X$ is in this interval, or $T^{kn_i+2n_i}_{kn_i+n_i} < T^{kn_i+2n_i}_{kn_i+4n_i}$, which necessarily occurs for $X$ at the site $(k+1)n_i$. Thus we have that the moment $\sigma_0$ at which it becomes certain that $\Gamma^{ln_i} ( e_i, n_i)^c$ occurs for some $l \in \left( \frac{n_{i+1}}{2n_i}-3, 4\frac{n_{i+1}}{n_i} \right)$ , either $B_1$ is at this moment impossible or we have that on every interval $(kn_i, (k+1)n_i)$ for $n_ik < X_\sigma-n_i$ and $k \geq \frac{n_{i+1}}{2n_i}$ there are at least $e_i \geq 1/2$ sites with $M_1$ cookies and for $x \geq X_{\sigma_0} + 4n_i$ the number of cookies is equal to $M$. We now define a stopping time subsequent to $\sigma_0$:
$$
\sigma_1 = \inf\{n > \sigma_0: X_n = X_{\sigma_0} + 9n_i \mbox{ or } X_n = X_{\sigma_0} - n_i^{5/4}\} \, .
$$

(We may assume that $n_i^{5/4} << n_{i+1}/2$.)
We have by Corollary \ref{cor2} and our choice of $M_1$ that for $n_1$ large
$$
\mathrm{P}(X_{\sigma_1} = X_\sigma + 9n_i) \ge 1- \epsilon^{\log_2(n_i^{1/4}/9)}
= 1 - (n_i^{1/4}/9)^{\log_2 \epsilon}
\ge 1 - \frac{1}{n_i^{-\frac{1}{5}\log_2 \epsilon}} \, .
$$
On the other hand we have on $\{X_{\sigma_1} = X_\sigma + 9n_i\} \cap\{ \mbox{There exists no } l, l^\prime \mbox{ with } \frac{n_{i+1}}{2n_i} \leq l, l^\prime \leq 4 \frac{n_{i+1}}{n_i}$ with $l^\prime - l > 3$ and $\Gamma^{n, l}(e_i,n_i)^c \cap \Gamma^{n, l^\prime} (e_i, n_i)^c\}$ that
$$
\Gamma^{X_{\sigma_0} + 8n_i} (e_i, n_i) \cap \Gamma^{X_{\sigma_0}+9n_i}(e_i, n_i) \ldots \cap \Gamma^{4n_{i+1}}(e_i, n_i)
$$
occurs. This implies
\begin{itemize}
\item[(i)] for every interval $[kn_i, kn_{i+1}]$ not intersecting $[X_{{\sigma_0}}-n_i^{5/4}, X_{{\sigma_0}}+9n_i)$ the number of sites with $M_i$ cookies will exceed $e_in_i$
\item[(ii)] upon hitting $X_{\sigma_0} + 9n_i$ the cookie walk passes from $X_{\sigma_0}+ kn_i$ to $X_{\sigma_0}+(k+2)n_i$ before returning to $X_{\sigma_0} + (k-1)n_i$ for all $k \geq 9$ and such that $k \leq 4\frac{n_{i+1}}{n_i} - \frac{X_{\sigma_0}}{n_i}$.
\end{itemize}
Thus we have on this event that
\begin{enumerate}
\item[(1)] After hitting $2n_{i+1}$ the cookie r.w. hits $4n_{i+1}$ before $n_{i+1}$
\item[(2)] upon hitting $4n_{i+1}$ the number of sites in $(n_{i+1}, 2n_{i+1})$ with $\geq M_0$ cookies is at least
$$
n_{i+1} e_i - 2n_i^{5/4} \ge n_{i+1} \left( e_i-\frac{3}{n_i^{1/4}} \right) \equiv n_{i+1} e_{i+1} \quad \textrm{(if } n_1 \textrm{ is (and therefore all } n_i \textrm{ are) large).}
$$
\end{enumerate}
Clearly if $e_1 = 3/4$ and $n_1$ is fixed large enough to ensure that $\sum^\infty_{i=1} \frac{3}{n_i^{1/4}} < 1/4$ then $\{e_i\}$ defined as above satisfy $e_i \geq 1/2$ for every $i$. The preceding show that the only way that event $\Gamma(e_{i+1}, n_{i+1})^c$ can occur is if either
\begin{enumerate}
\item[(i)] for two $l, l^\prime$ differing by more than 3 such that $l, l^\prime \in \left[0,4\frac{n_{i+1}}{n_i}\right]$ we have $\Gamma^{ln_i} (e_i, n_i)^c \cap \Gamma^{l^\prime_i}(e_i,n_i)^c$ 
\item[(ii)] $\{X_{\sigma_0} < 4n_{i+1}\} \cap \{X_{\sigma_1} =X_{\sigma_0}-n_i^{5/4}\}$, for $\sigma_0$ and $\sigma_1$ as above.
\end{enumerate}
This, together with (\ref{eq:twobad}), gives the inequality
$$
\mathrm{P}(\Gamma(e_{i+1}, n_{i+1})^c) \leq \left(4 \frac{n_{i+1}}{n_i} \mathrm{P}(\Gamma(e_i, n_i)^c) \right)^2
+ 4 \frac{n_{i+1}}{n_i} \frac{1}{n_i^{-\frac{1}{5}log_2 \epsilon}} \mathrm{P}(\Gamma(e_i, n_i)^c).
$$
By the induction hypothesis $\mathrm{P}(\Gamma(e_i, n_i)^c) \leq \frac{1}{n_i^{\zeta}}$ then
\begin{eqnarray*}
\mathrm{P}(\Gamma(e_{i+1}, n_{i+1})^c) &\le& \left(4 n_i^{1/2} \frac{C}{n_i^{\zeta}}\right)^2 + \frac{4 n_i^{1/2}}{n_i^{-\frac{1}{5}log_2 \epsilon}} \ \frac{C}{n_i^{\zeta}} \\
&\le& \frac{1}{n_{i+1}^{\zeta}} \left( \frac{16C^2}{n_i^{\zeta/2-1}} + \frac{4C}{n_i^{-\frac{1}{5}log_2 \epsilon - (\zeta+1)/2}} \right) \, .
\end{eqnarray*}
If $\epsilon$ is fixed sufficiently small and afterward $n_1$ is taken large, the expression under parenthesis in the rightmost term of the previous expression is bounded by 1 and then (\ref{eq:ind}) holds for $i+1$. $\square$

\bigskip \bigskip

\no \textbf{Counter-example} We wish to show through a counter example that the previous result cannot be generalized to ergodic environments, i.e., for a fixed $p>1/2$, we give an initial ergodic environment with mean number of cookies per site greater or equal to $M$ such that the ERW speed is $0$ independent on the choice of $M$.

Fix some $\frac{1}{2}<\epsilon<1$ and $0<\sigma<1$ such that $\sum_{n\ge 2} \sigma^n \ge M$. Now consider a sequence $(X_n)_{n\ge 1}$ of i.i.d random variables with the following distribution:
$$
X_n =
\left\{
\begin{array}{ll}
0 \ , & \textrm{with probability } \sigma, \textrm{ for } n\ge 2 \\
2^n \ , & \textrm{with probability } \gamma /4^{\epsilon n}, \textrm{ for } n\ge 2 ,
\end{array}
\right.
$$
where $\gamma$ is a normalization constant chosen so that $\gamma \sum _{n=2}^\infty \frac{1}{4^{n \epsilon}} \ + \ \sigma = 1$. With our choice for $\epsilon$, $E[X_n]$ is finite. Now let $(\mathcal{N}_x)_{x\ge 0}$ be the renewal counting process associated to $(X_n)_{n\ge 1}: N_x \ = \ \sup \{K: \sum _ {n=1}^k X_n \leq x\}$. We now replace $(\mathcal{N}_x)_{x\ge 0}$ by the corresponding translation invariant two sided process   $(\mathcal{N}_x)_{-\oo < x < +\oo}$ with $N_0 = 0$. We fix the initial environment as the following
$$
\omega(x,n) = \left\{
\begin{array}{ll}
\frac{3}{4} , & \textrm{for } 1 \le n \le (\mathcal{N}_{x+1} - \mathcal{N}_x ) \, , \\
\frac{1}{2} , & \textrm{otherwise} \, .
\end{array}
\right.
$$
The environment $\omega$ is stationary and ergodic, moreover, by the condition imposed on $\sigma$, the expected number of cookies per site is $E[\mathcal{N}_1] \ge M$. We have then that there exists a strictly positive constant $c$ (not depending on $n$) so that for all $n \geq 2$, the number of disjoint intervals which are initially empty of cookies and of length at least $2^n$ contained in $[0,M)$, $N_n(M)$, satisfies 
$$
N_n(M) \ \geq \ c 4^{-\epsilon n }M.
$$
Now the additional time $\tau_I$ for the cookie random walk to traverse such an interval, $I$, for the first time, after arriving at its leftmost endpoint is easily seen to satsify
$P( \tau_I > 4^n) > d >0$ for some constant not depending on $n$.  Thus we obtain by the law of large numbers
that
$$
\liminf_{M \rightarrow \infty} \frac{T_M}{M}  \geq c 4^{-\epsilon n }d 4^n.
$$
Since $\epsilon <1$ and $n$ is arbitrarily large, we have that the velocity must be zero.
%We also have that for $n_0$ sufficiently large
%\begin{eqnarray*}
%\mathrm{P} (\mathcal{N}_{2^n} \le 2^{n-1}) & = & \mathrm{P} \left( \sum_{i=1}^{2^{n-1}} X_i > 2^n \right) \ge 1 - \mathrm{P} (X_1 \le 2^n)^{2^{n-1}} \\
%& \ge & 1 - \left( 1 - \frac{\gamma}{3 \cdot 4^{\epsilon n}} \right)^{2^{n-1}} \ge \frac{C}{2^{\epsilon n}}
%\end{eqnarray*}
%for every $n\ge n_0$ and some constant $C>0$ not depending on $n$. Therefore we obtain as an application Lemma \ref{lem4} that
%\begin{eqnarray*}
%\lefteqn{\mathrm{P} ( \textrm{the exit time of } [-2^n,2^n] \textrm{ is greater or equal to } 2^{2n} ) } \\
%& & \ge E[ \textbf{1}\{\mathcal{N}_{2^n} \le 2^{n-1}\} \mathrm{P} ( \textrm{the exit time of } [-2^n,2^n] \textrm{ is greater or equal to } 2^{2n-1} | \mathcal{N}_t , -2^n \le t \le 2^t )  ] \\
%& & \ge \frac{c}{2^{\epsilon n}} \, .
%\end{eqnarray*}
%In this way we have that
%$$
%\frac{1}{v} = \lim_{n \ra +\oo} E\left[ \frac{T_{2^n}}{2^n} \right] \ge \lim_{n \ra +\oo} \frac{1}{2^n} \frac{2^{2n}}{2} \frac{c}{2^{\epsilon n}} = \lim_{n \ra +\oo} c 2^{n(1-\epsilon)} = +\oo \, ,
%$$
%and therefore $v=0$ as we have claimed.

\bigskip \bigskip

%%%%%%%%%%%%%%%%%%%%%%%%%%%%%%%%%%%%%%%%%%%%%%%%%%%%%%%%%%%%%%%
\section{Proof of (ii) in Theorem \ref{thm1}} 
\label{section:thm1(ii)}
%%%%%%%%%%%%%%%%%%%%%%%%%%%%%%%%%%%%%%%%%%%%%%%%%%%%%%%%%%%%%%%
\setcounter{equation}{0}

Consider the excited random walk $(X_n,\omega_n)_{n\in \N}$ starting at $(0,\omega^{M,p})$ for some fixed $p > 1/2$ and $M$ which satisfies $M(2p-1) \in (1,2)$. It is well known (see \cite{zeitouni} Lemma 2.1.17), that
\begin{equation}
\label{eq:T1}
\textrm{a.s}-\!\! \lim_{V \ra \oo} \frac{T_V}{V} = \left( \textrm{a.s}- \!\! \lim_{n \ra \oo} \frac{X_n}{n} \right)^{-1} = \frac{1}{\mu} ,
\end{equation}
which is the inverse of the walk's speed. By the dominated convergence theorem we have under truncation that
$$
\lim_{V \ra +\oo} \mathrm{E}_{(0,\omega^{M,p})} \left[ \frac{T_V}{V} \wedge m \right] = m \wedge \frac{1}{\mu} \, .
$$
However we are going to show that for $m$ large
$$
\liminf_{V \ra +\oo} \mathrm{E}_{(0,\omega^{M,p})} \left[ \frac{T_V}{V} \wedge m \right] \ge C m^\gamma \, ,
$$
for some $0<\gamma<1$, which implies that $\mu = 0$.

\smallskip

First we explain how the proof works. Accordingly we fix $v,v^\prime$ so that $1< M(2p-1) < 1+v < 1 + v^\prime < 2$. Our first aim will be to show that with sufficiently large probability for a fixed large density of $x$'s
\begin{equation}
\label{eq:T2}
\mathrm{P}_{(0,\omega^{M,p})} [ T_{x+1} - T_x \ge 2^{2r} | \mathcal{F}_{T_x} ] \ge c 2^{-r(1+ v^\prime)} 
\end{equation}
for some $c>0$ and every $r$ sufficiently large. For $x$ satisfying (\ref{eq:T2}) and a sufficiently large $m$ and $r$ such that $2^{2r} \le m \le 2^{2(r+1)}$ 
$$
\mathrm{P}_{(0,\omega^{M,p})} [ T_{x+1} - T_x \ge m | \mathcal{F}_{T_x} ] \ge \frac{c}{m^\alpha}
$$
for $\alpha = (1+ v^\prime)/2 < 1$. From this last inequality we will have that the limit in \ref{eq:T1} is $\infty$ and so the ERW velocity will be zero. 

\smallskip

The proof of (\ref{eq:T2}) is based on applications of Lemma \ref{lem5} and in conformity with its statement we fix arbitrarily $\epsilon$ small and $N \le (\epsilon (2p-1))^{-1/2}$ so that in notation of Lemma \ref{lem5} 
$$
2^{-(1+v)}(1+c_N) > 2^{(1+v^\prime)} \, . 
$$
and let $K$ be sufficiently large so that 
\begin{equation}
\label{eq:T3}
\max_{k\ge K} C_k \geq 2^{-(1+v^\prime)} \, ,
\end{equation}
We start by stating some direct implications of Lemma 10 in \cite{zerner} which gives that the sequence
$$
((X_{\tau_{k, m}} -k)_{m \geq 0})_{k \geq 0} = (Z_k)_{k \geq 0}
$$
is stationary ergodic in $k$, where $(\tau_{k,m})_{m \ge 1}$ is defined by $\tau_{k,0} = -1$ and $\tau_{k,m+1} = \inf\{ n > \tau_{k,m} : X_n \ge k \}$. We can without loss of generality assume that $(Z_k)_{k \geq 0}$ extends to an ergodic two sided process $(Z_k)_{-\infty < k < \infty}$. The implications we mentioned are:

\medskip

\no \textbf{Fact 1:} Consider the functional $f(Z_k)$ which is the number of cookies left at site $k$ finally. We have \cite{zerner}, Lemma 2 that $a.s$ for every $k$, 
\begin{equation}
\label{eq:mconv1}
\lim_{n\ra +\oo} \frac{1}{n+1} \sum^k_{x = -n+k} f(Z_x) = \frac{1}{(2p-1)} ((2p-1)M-1) < \frac{v}{2p-1} \, .
\end{equation}
This implies that a.s. for every $x$, there exists $L_x < \ \infty$ such that for every $r \ge L_x$, 
\begin{eqnarray*}
\sum^{x-\frac{2^{r}(i-1)}{N}-2^{r}}_{y = x - \frac{2^{r}i}{N}-2^{r}}  f(Z_y) \leq \frac{2^{r}}{N} \frac{v}{2p-1}
\end{eqnarray*}
for every $i = 1, 2, \ldots,  N$. Furthermore we can also fix $L$ such that a.s.
$$
\textrm{a.s-}\!\!\!\!\! \lim_{n\ra +\oo} \frac{1}{n} \sum^n_{y=0} I_{L_y \leq L} > 9/10.
$$

\medskip

\no \textbf{Fact 2:} Let $h(Z_x, l)$ be equal to $1$ if the number of excursions from $x$ to $(-\infty, x)$ is equal to $l$, otherwise its value is $0$. For every $l$ there exists $c(l) < 1$ such that  
\begin{equation}
\label{eq:mconv2}
\textrm{a.s-}\!\!\!\!\! \lim_{n \ra +\oo} \frac{1}{n+1} \sum^k_{x=-n+k} h (Z_x, l) = c(l).
\end{equation}
Note that since $M(2p-1)>1$ implies that X is transient by \cite{zerner}, Theorem 12, then almost surely there exists $l$ such that $h(Z_x, l)=1$ and thus $\sum_l c(l) = 1$. Moreover we can pick $k_0$ so large that the density of $x$ making at most $k_0$ excursions from $x$ to $(-\infty, x)$ is at least 9/10, i.e.
$$
\textrm{a.s-}\!\!\!\!\! \lim_{N \ra +\oo} \frac{1}{N} \sum^N_{x=0} \sum^{k_0}_{l=0} h(Z_x, l) \geq 9/10.
$$

\medskip

\no \textbf{Fact 3:} We consider for $x$ fixed the indicator function $g_R(x)$ of the event that among the first $k_0$ excursions of $X$ from $x$ to $(-\infty, -x)$ the value $x-2^R$ is realized. By the Markov property and simple majoration of $X$ by appropriate random walks we can fix ${L^\prime}>0$ sufficiently large such that for every $x$ we have $E \left[ g_{L^\prime}(x) \right] < 100^{-1}$ and so with probability greater than $\frac{9}{10}$, for $n$ large,
$$
\frac{1}{n} \sum^{n-1}_0 g_{L^\prime}(x) \leq \frac{1}{10}.
$$

\medskip

What we get from Facts 1-3 is that for $n$ large with probability at least $\frac{9}{10}$, a density of $\frac{7}{10}$ of sites $x$ are going to satisfy: $L_x \le L$ and $T^x_{x-2^{L^\prime}}=+\oo$, i.e., at the time $x$ is visited for the first time the ERW satisfies  
$$
\sum^{x-\frac{2^{r}(i-1)}{N}-2^{r}}_{y = x - \frac{2^{r}i}{N}-2^{r}}  f(Z_y) \leq \frac{2^{r}}{N} \frac{v}{2p-1}
$$
for every $r \ge L$ and does not visit $x - 2^{L^\prime}$ again. We say that the environment is \emph{good} on the right of $x - 2^{R}$ for $R=max\{L,L^\prime\}$, if there is at most $\frac{2^R}{N} \cdot \frac{v}{2p-1}$ cookies on each interval 
$$
\left( x - \frac{2^l i}{N} - 2^l , x - \frac{2^l (i-1)}{N} - 2^l \right) \, .
$$
for $i=1,...,N$ and $l \ge R$. We say that a site $x$ is \emph{good} if the environment on the right of $x - 2^{R}$ is good.

\smallskip

Now fix $\beta$ as the infimum, over all possible environments with at most $M$ cookies per site, of the probability that the ERW starting from 0 hits $- 2^{R}$ before $1$ eating all the cookies in the interval $(-2^{R},0]$ (assume without loss of generality that $R > 0$). Suppose that $x$ is good, then by the Markov property for $(X_n,w_n)_{n \ge 1}$ and Lemma \ref{lem4} (with $\gamma$ taken as $\nu/(2p-1)$) we have that
$$
\mathrm{P}_{(0,\omega^{M,p})} [ T_{x+1} - T_x \ge 2^{2r} | \mathcal{F}_{T_x} ] \ge \beta \, c_0 \, \mathrm{P} \left[ \left. T_{x+1} \le T^{x-2^{R}}_{x-2^{r+1}} \right| \mathcal{F}_{T_x} \right] \, ,
$$
for some universal constant $c_0$ not depending on $r$ and $x$. The initial environment, $\omega^{M,p}$, conditioned to $\mathcal{F}_{T_x}$, for a good $x$, translated by $2^{R}-x$ satisfies the hypotheses of Lemma \ref{lem5} with $b=v/(2p-1)$, and thus we can apply the lemma recursively to the events obtained by translating $\Omega$ by $x-2^{R}$, $x-2^{R+1}$, ... , $x-2^{1+r}$ to obtain that 
$$
\mathrm{P}_{(0,\omega^{M,p})} \left[ \left. T_{x+1} \le T^{x-2^{R}}_{x-2^{r}} \right| \mathcal{F}_{T_x} \right] \ge \prod_{k=R}^{r+1} C_k \, .
$$
We may, without loss of generality, suppose $R \ge K$ from (\ref{eq:T3}) and thus
$$
\mathrm{P}_{(0,\omega^{M,p})} [ T_{x+1} - T_x \ge 2^{2r} | \mathcal{F}_{T_x} ] \ge \beta \, c \, \prod_{k=R}^r C_k \ge c^\prime 2^{-r(1+v^\prime)} \, .
$$
Therefore, we have now that with probability at least $\frac{9}{10}$ for large $V$
\begin{equation}
\label{eq:Tcond}
\mathrm{P}_{(0,\omega^{M,p})} [ T_{x+1} - T_x \ge 2^{2r}| \mathcal{F}_{T_x} ] \ge  C \, 2^{-r(1+v^\prime)}
\end{equation}
for $\frac{7}{10}$ of $x$ in $[0,V]$. Hence we write
$$
\mathrm{E}_{(0,\omega^{M,p})} \left[ \frac{T_V}{V} \wedge m \right] = \mathrm{E}_{(0,\omega^{M,p})} \left[ \frac{1}{V} \sum_{x=0}^{V-1} (T_{x+1}-T_x) \wedge m \right]
$$
which, by (\ref{eq:Tcond}), for $V$ large is bounded below by 
\begin{eqnarray}
\lefteqn{\!\!\!\!\!\!\!\!\!\!\!\!\!\!\!\!\!\!\!\!\!\!\!\!\!\!\!\!\!\!\!\!\!\!\!\!
\mathrm{E}_{(0,\omega^{M,p})} \left[ \frac{1}{V} \sum_{i=1}^{V-1}  m \, \mathrm{P}_{(0,\omega^{M,p})} [ T_{x+1} - T_{x} \ge m | \mathcal{F}_{T_{x}}]  \right] \ge} \\
& & \ge \frac{C}{V} \left(\frac{7}{10}(V-2^{2R})-m\right) m^{1- \frac{1+v^\prime}{2}} \, ,
\end{eqnarray}
and then, for $\gamma = 1- \frac{1+v^\prime}{2} > 0$, we have that
$$
\mathrm{E}_{(0,\omega^{M,p})} \left[ \frac{T_V}{V} \wedge m \right] \ge \frac{C}{V} \left(\frac{7}{10}(V-2^{2R})\right) m^\gamma \ge C^{\prime \prime} m^\gamma \, ,
$$
for $V$ large. $\square$

\newpage

%%%%%%%%%%%%%%%%%%%%%%%%%%%%%%%%%%%%%%%%%%%%%%%%%%%%%%%%%%%%%%%%%%%%%%%%%%%%%%%%%
\section{Proof of Theorem \ref{thm2}}
\label{section:thm2}
%%%%%%%%%%%%%%%%%%%%%%%%%%%%%%%%%%%%%%%%%%%%%%%%%%%%%%%%%%%%%%%%%%%%%%%%%%%%%%%%%

\subsection{Preliminaries}
Let $\epsilon,\kappa>0$ let $W$ be a nearest neighbour random walker on $\mathbb{Z}$ which jumps to the right with probability $1/2 + \epsilon$, when it is at position $x\in [-\kappa/\epsilon,\kappa/\epsilon]$, while for $x\not\in [-\kappa/\epsilon,\kappa/\epsilon]$ it moves like a symmetric random walker. Let $L,v>0$ be positive integers and consider the following events:
\[
A_1=A_1(\epsilon,L,\kappa)=\Big\{ W\mbox{ touches $\frac{L+\kappa}{\epsilon}$ for the first time before it touches $\frac{-L+\kappa}{\epsilon}$}\Big\}\,;
\]
\[
A_2=A_2(\epsilon,L,\kappa,v)=
\]
\[
\Big\{ \exists\,x\in[\frac{-\kappa}{\epsilon},\frac{\kappa}{\epsilon}]\,:\, W\mbox{ has visited $x$ more than $\frac{v}{\epsilon}$ times before it hits $\frac{L+\kappa}{\epsilon}$}\Big\}\,.
\]

\begin{lemma}\label{comp0}
There exist constants $L_0,\epsilon_0,\kappa_0>0$ and $c_0>0$ such that for all $\epsilon<\epsilon_0$ and $\kappa<\kappa_0$ so that for event $A:=A_1$ defined by the parameters $L_0,\epsilon,\kappa$, we have that for all $x\in [-\kappa/\epsilon,\kappa/\epsilon]$:
\[
P^{x,\bar{\omega}}(A)> \frac{1}{2}+c_0\kappa\,.
\]
\end{lemma}
This can be seen by exlicit calculation for what is a birth and death process (see e.g. \cite{HPS}).  We have also that for this birth and death process, the excursion theory is easily analysed and we obtain via a simple chaining argument
\begin{lemma}\label{comp2}
For the $L_0$ fixed in Lemma \ref{comp0}, there exists $v < \infty$ so that for all $\epsilon<\epsilon_0$ and $\kappa<\kappa_0$ (the constants of Lemma \ref{comp0}) the  event $A:=A_2$ defined by the parameters $L_0,\epsilon,\kappa,v$, we have that for all $x\in [-\kappa/\epsilon,\kappa/\epsilon]$:
\[
P^{x,\bar{\omega}}(A^c)< e^{-\frac{1}{\kappa}}.
\]
\end{lemma}
Putting the two lemmas together we find (possibly at the price of reducing further $\kappa _0$)
\begin{lemma}\label{comp1}
There exist constants $L_0,\epsilon_0,\kappa_0>0$ and $c_1>0$ such that for all $\epsilon<\epsilon_0$ and $\kappa<\kappa_0$ so that for event $A:=A_1\cap A_2^c$ defined by the parameters $L_0,\epsilon,\kappa$, we have that for all $x\in [-\kappa/\epsilon,\kappa/\epsilon]$:
\[
P^{x,\bar{\omega}}(A)> \frac{1}{2}+c_1\kappa\,.
\]
\end{lemma}

\begin{remark}  The object is to apply the above to a cookie random walk.  
The result above immediately implies that for a cookie random walk (with bias parameter $\frac12 + \epsilon $), at least $\frac{v}{\epsilon}$ cookies at each point in interval
$[-\frac{\kappa}{\epsilon},\frac{\kappa}{\epsilon}]$ and with no cookies outside of this interval, we have that
$P^{x,\bar{\omega}}(A)> \frac{1}{2}+c_1\kappa $.
 Since the existence of cookies to the left of $\frac{-\kappa}{\epsilon}$ can only help event $A_1$ while not affecting at all event $A_2$, these have no effect on the above bounds.  Likewise the existence of cookies to the right of $\frac{\kappa}{\epsilon}$ can only increase the probability of event $A$ and so Lemma \ref{comp1}
 applies to cookie random walks in some generality.
 \end{remark}
\subsection{Coupling}
By Lemma \ref{comp1}, one can fix $L_0>0$ sufficiently large and $\epsilon_0, \kappa>0$ sufficiently small so that for all $\epsilon < \epsilon_0$ we have 
\[
P^{x,\omega}(A)\geq 1/2 + c_1 \kappa\mbox{ for all }x\in [-\kappa/\epsilon,\kappa/\epsilon]\,.
\]
From Theorem \ref{thm1}, there exists $M_0=M_0(\kappa)>1$ so that ERW with $M_0$ cookies and bias $p=1/2 + c_1\kappa$ has a strictly positive drift. From now on we denote by $X$ the ERW with these parameters. Let $X^\epsilon$ be a walker with environment given by $\omega^{M_0 v/\epsilon,1/2+\epsilon}$. With Lemma \ref{comp1} in hand we show that:
\begin{proposition}\label{coup}
There exists a coupling between $X$ and $X^\epsilon$ with the following properties: there exist an increasing sequence of stopping times $(\tau_n^{\epsilon})_{n\geq 1}$ with 
\begin{equation}\label{coup1}
\limsup_n \frac{\tau_n^\epsilon}{n} < \infty\,
\end{equation}
and a random function $l:\mathbf{N}^+\to\mathbf{N}^+$ with $l(n+1)-l(n)\geq 1$ such that
\begin{equation}\label{coup2}
 \frac{\epsilon X^\epsilon_{\tau_n^\epsilon}}{L_0}-X_{l(n)}\geq -\kappa_1\,.
\end{equation}
\end{proposition}
We first show how this result leads to Theorem \ref{thm2}.\\
\no \textbf{Proof of Theorem \ref{thm2}} 
First note that to prove the desired result it is only necessary to treat $ \epsilon$ small, so restricting to $\epsilon < \epsilon_0$ of Lemma \ref{comp1} is certainly legitimate.  By \ref{coup1} in Proposition \ref{coup} (and that $l(n)\geq n$),
\[
\frac{\epsilon}{L}\liminf_n\frac{X^\epsilon_{\tau_n^\epsilon}}{\tau_n^\epsilon}\geq \liminf_n\frac{X_{l(n)}}{\tau_n^\epsilon}\geq \liminf_n\frac{X_{l(n)}}{l(n)}\frac{n}{\tau_n^\epsilon}\,.
\]
Since $X$ has a strictly positive speed, together with \ref{coup2} this yields that $X^\epsilon$ has a strictly positive speed. $\square$

\bigskip

\no \textbf{Proof of Proposition \ref{coup}} Consider the collection of intervals  $\{rL/\epsilon +[-\kappa/\epsilon,\kappa/\epsilon]\,:\,r\in\mathbb{Z}\}$. Let $\tau_0^\epsilon=0$ and $I(0)=[-\kappa/\epsilon,\kappa/\epsilon]$ and assume that $\tau_0^\epsilon,\dots,\tau_n^{\epsilon}$ were defined. These stopping times will be chosen so that for each $n\geq 0$ there will be a unique interval  $I(n)\in\{rL+[-\kappa/\epsilon,\kappa/\epsilon]\,:\,r\in\mathbb{Z}\}$ so that $X^\epsilon_{\tau^\epsilon_n}\in I(n)$.  We also introduce the process $(Z_n\,:\,n\geq 0)$ by setting $Z_n=r$ if $I(n)=rL/\epsilon+[-\kappa/\epsilon,\kappa/\epsilon]$. Now, define $\tau_{n+1}^\epsilon$ as follows: If $\sum_{i=0}^{n}I_{Z_i=Z_n}\leq M_0 -1$ then let $\tau^\epsilon_{n+1}$ be the first time after $\tau_n^\epsilon$ that
\begin{itemize}
\item there exists $x\in I(n)$ so that the number of visits to $x$ is greater than $v/\epsilon$ ,  
\item or $X^\epsilon$ touches $(L(Z_n-1)+\kappa)/\epsilon$ ,
\item or $X^\epsilon$ touches $(L(Z_n+1)-\kappa)/\epsilon$ if $\sum_{i}^{n}I_{Z_i=Z_n+1}\leq M_0 -1$ ,
\item or $X^\epsilon$ touches $(L(Z_n+1)+\kappa)/\epsilon$ if  $\sum_{i}^{n}I_{Z_i=Z_n+1}\geq M_0$ ;
\end{itemize}
If $\sum_{i=0}^{n}I_{Z_i=Z_n}\geq M_0$ then let $\tau^\epsilon_{n+1}$ be the first time after $\tau_n^\epsilon$ that
\begin{itemize}
\item $X^\epsilon$ touches $(L(Z_n-1)+\kappa)/\epsilon$ ,
\item or $X^\epsilon$ touches $(L(Z_n+1)-\kappa)/\epsilon$ if $\sum_{i}^{M_0}I_{Z_i=Z_n+1}\leq M_0 - 1$ ,
\item or $X^\epsilon$ touches $(L(Z_n+1)+\kappa)/\epsilon$ if  $\sum_{i}^{M_0}I_{Z_i=Z_n+1}\geq M_0$ .
\end{itemize}

From this definition it should be clear that
\[
\limsup_n \frac{\tau^\epsilon_n}{n}< K_\epsilon\,
\]
(where $K_\epsilon$ is a finite constant depend on $\epsilon$). Notice that if $\sum_{i=0}^{n}I_{Z_i=Z_n}=c\leq M_0 - 1$ then every site in $I(n)$ has at least $v(M_0 - c)/\epsilon > v/\epsilon$ cookies. We also make the following observations:
\begin{enumerate}
\item If $Z_{n+1}=Z_n-1$ then $X^{\epsilon}_{\tau_{n+1}^\epsilon}= (Z_n L+\kappa)/\epsilon$;
\item If $\sum_i^{n} I_{Z_i = Z_n} \leq M_0 $ then, by Lemma \ref{comp1} (and the remark following it), the probability that $Z_{n+1}=Z_n+1$ (or equivalently $X^\epsilon_{\tau_{n+1}^\epsilon}\in (I(n)+L/\epsilon)$) is greater than $1/2 +c_1\kappa$;
\item If $\sum_{i=0}^{n}I_{Z_i = Z_n} > M_0 $ then, by the rules for $\tau_n^\epsilon$,  $X^\epsilon_{\tau_{n+1}^\epsilon}= (Z_n L+\kappa)/\epsilon$ and so the conditional probability given $\mathcal{F}_{\tau_n^\epsilon}$ that $Z_{n+1}=Z_n + 1$ is at least $1/2$  (here $\mathcal{F}_.$ refers to the filtration for the process $X^\epsilon_.$).
\end{enumerate}

It is important to realize the purpose of the different rules when the number of previous visits to the current site goes from $M_0 -1$ to $M_0$:  it would be possible that when $Z_n=r$ for the $M_0$ time that $X_{\tau_n^\epsilon}^\epsilon$ were to the left of $Z_n L/\epsilon$. This would mean that that the conditional probability of advancing might be less than $1/2$. It is to make the "bridge" for the subsequent regime where (3) holds that we change the rules for the stopping times.

We now introduce our comparison process: the idea is to consider $X=(X_n\,:\,n\geq 0)$ a $M_0$ cookie for bias  $p=1/2+c_1\kappa$. We do not couple so that $X_n\leq Z_n$ for all $n$ but rather we will construct $l(n)$ so that  $l(0)=0$, $l(n+1)-l(n)\geq 1$ and $X_{l(n)}=Z_n$ quite simply. We will also have that for all $r\in \mathbb{Z}$
\[
\sum_{i=0}^{l(n)} I_{X_i=r}\geq \sum_{i=0}^{n} I_{Z_i=r}\,.
\]
Thus, given $l(n)$, we have either 
\[
\sum_{i=0}^{n}I_{Z_i=Z_n} > M_0 \mbox{ or } \sum_{i=0}^{n}I_{Z_i=Z_n}\leq M_0\,.
\]
In the former case, by our inductive hypothesis,
\[
P(Z_{n+1}=Z_n + 1\mid \mathcal{F}_{\tau_n^\epsilon})\geq \frac{1}{2} =P(X_{l(n+1)}=X_{l(n)} +1\mid \sigma(X_0,\dots,X_{l(n)}))
\]
and so we can couple them so that
\[
\big\{X_{l(n+1)}=X_{l(n)} + 1\big\}\subseteq \big\{Z_{n+1}=Z_n + 1\big\}\,.
\]
In the latter case we have (by item (2)) 
\[
P(Z_{n+1}=Z_n + 1\mid \mathcal{F}_{\tau_n^\epsilon})\geq \frac{1}{2} + c_1\kappa \geq P(X_{l(n+1)}=X_{l(n)} +1\mid \sigma(X_0,\dots,X_{l(n)}))
\]
and again we can couple them.

To define $l(n)$ we simple set $l(n+1)=l(n) +1$ if $Z_{n+1}=X_{n+1}$. If not, we take $l(n+1)=\inf\{ k\geq l(n)\,:\, X_k= Z_{n+1}\}$. Thus $X_{l(n)}=Z_n$ and $\sum_{i=0}^{l(n)} I_{X_i=r}\geq \sum_{i=0}^{n} I_{Z_i=r}$ for all $r$, as we required. Together with the definition of $l(n)$, item (1) implies \ref{coup2} and the proof of Proposition \ref{coup} is complete. $\square$

\bibliographystyle{alpha}

\begin{thebibliography}{99}

\bibitem{kozma} Kozma G.: Excited Random Walks in three dimensions has positive speed, 

\bibitem{HPS} Hoel, P.G., Port, S. and Stone C.: Introduction to Stochastic Processes, Waveland Pr. Inc.

\bibitem{benjamini-wilson} Benjamini I., Wilson D.: 
Excited Random walk, Elect. Comm. Probab. 8, 86-92 (2003) 

\bibitem{zeitouni} Notes on Saint Flour Lectures 2001, preprint, \\
http://www.ee.technion.ac.il/~zeitouni/ps/notes1.ps.

\bibitem{zerner} Zerner M.: 
Multi-excited random walks on the integers

\end{thebibliography}

\end{document}